\documentclass{article}
\usepackage{graphicx}
\usepackage{cite}
\usepackage{mathtools}
\usepackage{xcolor}
\usepackage{siunitx}
\usepackage{tikz-cd}
\usepackage{microtype}
\usepackage{amsmath}
\usepackage{amssymb}
\usepackage{amsthm}

\usepackage{authblk}
\usepackage{hyperref}
\usepackage[a4paper, top=3cm, bottom=3cm, left=3cm, right=3cm]{geometry}

\title{The Induced Matching Distance: A Novel Topological Metric with Applications in Robotics
\footnote{Partially supported by
REXASI-PRO H-EU
project, call HORIZON-CL4-2021-HUMAN-01-01, Grant agreement ID: 101070028.}
\footnote{Álvaro Torras-Casas contract is funded by the French Agence Nationale de la Recherche through the project
reference ANR-22-CPJ1-0047-01}
\footnote{This work has been submitted to the \href{https://kam.mff.cuni.cz/conferences/eurocg2025/}{EuroCG 2025 Conference}.}
}

\author[1]{Javier Perera-Lago}
\author[2]{Álvaro Torras-Casas}
\author[3]{Jérôme Guzzi}
\author[1]{Rocio Gonzalez-Diaz}

\affil[1]{Universidad de Sevilla, Seville, Spain,  \texttt{$\{$jperera,rogodi$\}$@us.es}
}

\affil[2]{Inserm, INRAE, Centre for Research in
 Epidemiology and Statistics (CRESS),
 Université Paris Cité and Université
 Sorbonne Paris Nord
\texttt{alvaro.torras-casas@inserm.fr}
}

\affil[3]{SUPSI, IDSIA, Lugano, Switzerland,
\texttt{jerome.guzzi@idsia.ch }
}

\date{}

\DeclareMathOperator{\PH}{PH}
\DeclareMathOperator{\VR}{VR}
\DeclareMathOperator{\Ho}{H}
\DeclareMathOperator{\B}{B}
\DeclareMathOperator{\Rep}{Rep}
\DeclareMathOperator{\triplet}{triplet}
\DeclareMathOperator{\TMT}{TMT}

\newcommand{\R}{\mathbb{R}^+}
\newcommand{\N}{\mathbb{N}}
\newcommand{\Z}{\mathbb{Z}}
\newcommand{\cM}{\mathcal{M}}
\newcommand{\scst}{\scriptscriptstyle}
\newcommand{\rhoZ}{\rho^{\scst Z}}

\newcommand{\cMXZ}{\cM_f^0}
\newcommand{\veps}{\varepsilon}

\newcommand{\dM}{d_{f_0}}
\newcommand{\dX}{d_{\scst X}}
\newcommand{\dZ}{d_{\scst Z}}
\newcommand{\SZ}{S^{\scst Z}}
\newcommand{\SX}{S^{\scst X}}
\newcommand{\mXdelta}{m^{\scst X^{\delta}}}
\newcommand{\mZ}{m^{\scst Z}}
\newcommand{\mX}{m^{\scst X}}
\newcommand{\fbullet}{f_{\bullet}}
\newcommand{\ka}{\kappa^a}
\newcommand{\kb}{\kappa^b}
\newcommand{\kc}{\kappa^c}
\newcommand{\kinfty}{\kappa^\infty}
\newcommand{\Xdelta}{X^{\delta}}

\newcommand{\dXp}{d_{\scst X'}}

\newtheorem{theorem}{Theorem}[section]
\newtheorem{lemma}[theorem]{Lemma}

\begin{document}

\maketitle

\begin{abstract}
This paper introduces the induced matching distance, 
a novel topological metric designed to compare discrete structures represented by a symmetric non-negative function. We apply this notion to analyze agent trajectories over time. We use dynamic time warping to measure trajectory similarity and compute the 0-dimensional persistent homology to identify relevant connected components, which, in our context, correspond to groups of similar trajectories. To track the evolution of these components across time, we compute induced matching distances, which preserve the coherence of their dynamic behavior. We then obtain a 1-dimensional signal that quantifies the consistency of trajectory groups over time. Our experiments demonstrate that our approach effectively differentiates between various agent behaviors, highlighting its potential as a robust tool for topological analysis in robotics and related fields.
\end{abstract}

\section{Persistence Matching Distance
induced by Bijections}
From this point forward, we consider sets of $n$ points 
 $Z=\{z_1,z_2,\dots,z_n\}$ 
 together with a 
 symmetric non-negative function
 $\dZ\colon Z\times Z\rightarrow \R$,
where  $\R$ is the set of non-negative real numbers.
 Here, we remark that $\dZ$ might admit zero values between different points and does not need to satisfy the triangle inequality; as this is
 the case with the dynamic time warping measure considered in Section~\ref{sec:DTW}.
Since we are only interested in the 0-dimensional persistent homology, we work with the 1-skeleton  of the Vietoris-Rips filtration of $Z$ and we fix 
 $\Z_2$ as the ground field.

 The \emph{1-skeleton  of the Vietoris-Rips filtration} of $Z$,
denoted by $\VR(Z)$, 
is a family of graphs $\{\VR_r(Z)\}_{r\in\R}$ whose vertex set is $Z$. We set $\VR_0(Z)$ to have no edges, and, for all $r>0$, each graph $\VR_r(Z)$ contains the edges $[z,z']$ such that $\dZ(z,z')\leq r$.

Now, given $r\in \R$,
we denote by $\pi_0(\VR_r(Z))$ the quotient space $Z/\sim_r$, where two points $z,z' \in Z$ are such that $z\sim_r z'$ if and only if both $z$ and $z'$ lie in the same connected component from $\VR_r(Z)$.
We write $[z]$ to refer to the class of $z$ in $Z/\sim_r$, where $r \in \R$ is implied by the context.
The \emph{0-homology group}
$\Ho_0(\VR_r(Z))$
is the free $\Z_2$-vector space generated by the set $\pi_0(\VR_r(Z))$. 
The \emph{0-persistent homology} of $\VR(Z)$, denoted as $\PH_0(Z)$,
is given by the set 
of 0-homology groups 
$\big\{
\Ho_0(\VR_r(Z))
\big\}_{r\in \R}$ 
and the set of linear maps, 
$
\big\{
\rhoZ_{rs}\colon \Ho_0(\VR_r(Z)) \to\Ho_0(\VR_s(Z))
\big\}_{r\leq s}
$
that are induced by the inclusions
$\VR_r(Z)\subseteq \VR_s(Z)$
for all $r\leq s$.

The 0-persistent homology is an example of \emph{persistence module}. Formally, a persistence module $V$ is a set of free finite $\Z_2$-vector spaces $\{V_r\}_{r\in \R}$ together with a set of linear maps, called \emph{structure maps} $\{\rho_{rs}^{\scst V}\colon V_r\rightarrow V_s\}_{r\leq s}$.
Another example of persistence modules is
the \emph{interval module}, $\ka$, for 
$a>0$,
where $\ka_r= \Z_2$ for all $r <a$ and $\ka_r= 0$ otherwise, with structure maps 
$\{\rho^{\ka}_{rs}\colon \ka_r\to \ka_s\}_{r\leq s}$ being the identity map whenever $s\leq a$ and  the zero map otherwise.
In addition, we define $\kappa^0$ to be such that $\kappa^0_0= \Z_2$ and $\kappa^0_r= 0$ for all $r>0$.
Also, we denote by  $\kinfty$  
the persistence module that consists of 
$\kinfty_r= \Z_2$ for all $r\in R$; with structure maps $\{\rho^{\kinfty}_{rs}\colon\kinfty_r\to \kinfty_s\}_{r\leq s}$ being the identity maps.
Finally, a \emph{persistence morphism} $f\colon V\rightarrow U$ between persistence modules $V$ and $U$ is a set of linear maps $\big\{f_r\colon V_r\rightarrow U_r\big\}_{r \in \R}$ that commute with the structure maps 
$\rho^{\scst V}$ of $V$ and 
$\rho^{\scst U}$ of $U$. 
f $f_r$ is an isomorphism for all $r \in \R$ then $f$ is called a \emph{persistence isomorphism}.
An example of a persistence morphism is
$\kappa^{a} \rightarrow \kappa^{b}$, for $a\geq b$, where $\kappa^{a}_r \rightarrow \kappa^{b}_r$
consists of the identity map whenever $b\leq r\leq a$ and the zero map otherwise. An example of persistence isomorphism is
$\kinfty\rightarrow \kinfty$
where $\kinfty_r \rightarrow \kinfty_r$
consists of the identity map for all $r\in\R$.

\subsection{Connected Components, Barcodes and Triplet Merge
Trees}

Intuitively, $\PH_0(Z)$ encapsulates the evolution of connected components in $\VR(Z)$. 
This way,  all classes in $\PH_0(Z)$ are \emph{born} at $0$
since
$
\pi_0(\VR_0(Z))=
\big\{ [z_1],[z_2],\dots,[z_n] \big\}
$.
As the filtration parameter increases, $\PH_0(Z)$ records the death values of such classes.
Specifically, a class $[z_j]\in \pi_0(\VR_0(Z))$ is said to \textit{die at value $b>0$} 
if: 
\begin{itemize}
    \item[1)] $
\rhoZ_{0\ell}([z_j])=[z_j]$ for all $\ell\in\R$ with $\ell<b$.
\item[2)] 
$
\rhoZ_{0b}([z_j])=[z_i]
$, 
for some $i<j$.
That is, 
$[z_j]+[z_{i}]\in 
\ker\rhoZ_{0b}$.
\end{itemize}
In addition, we say that $[z_j]\in \pi_0(\VR_0(Z))$ \emph{dies at 
$0$} whenever
there exists $i<j$ such that 
$\dZ(z_i,z_j)=0$. 
Finally, observe that
the component $[z_1]$ never dies.
\\
A handy way to track the evolution of connected components is via 
\emph{triplet merge trees}
~\cite{triplets}: 
  \[\mbox{$\TMT(Z)=\big\{\triplet(j)\colon j\in\{2,3,\dots,n\}\big\} \subset Z\times \R \times Z\,\big\}$,}\]
where $\triplet(j)=(z_j, b_j, z_i)$ is
such that $[z_j]\in \pi_0(\VR_0(Z))$ dies at value $b_j\geq0$ and $\rhoZ_{0b_j}([z_j]) = [z_i]$, where $i\in\{1,2,\dots,j-1\}$ is as small as possible.
Using $\TMT(Z)$, we relate the generators of $\Ho_0(\VR_0(Z))$ to death values of connected components.

To store such death values, we use \emph{barcodes}.
The barcode $\B(Z)=(\SZ,\mZ)$, 
is a multiset where 
 $\SZ$ is the set of death values $b>0$, and $\mZ(b)$ is the number of generators of $\Ho_0(\VR_0(Z))$ that die at $b$, for a fixed $b>0$.
 The \emph{representation of}  $\B(Z)$ is
 $
\Rep\, \B(Z) = \big\{ \; (b,\ell) \;\vert\; b \in \SZ
\mbox{ and } 
\ell\in \{1,2,\dots,\mZ(b)\}\; \big\}
$.  
The following is a well-known fact
(see Theorem~1.2 of~\cite{decomposition2})
that we have adapted to our case. Its proof can be consulted in Appendix~\ref{prooflemma:decomposition-modules}.

\begin{lemma}~\label{lemma:decomposition-modules}
   There is a persistence isomorphism
    \[
       \begin{tikzcd}
       \PH_0(Z) 
    \ar[r, "\simeq"] & 
    \left(\oplus_{b \in \R} \oplus_{\ell \in   
    \{1,2,\dots,\mZ(b)\}} \kb\right) 
    \oplus 
    \kinfty.
    \end{tikzcd} 
   \]
\end{lemma}

\subsection{Induced Block Functions
and non-Expansive Maps}\label{sec:non-expansive}

To introduce  the notion of block functions that will be used later to define induced matching distances, we need the 
 operators $\ker^\pm_{b}$  defined as follows. For all $b>0$,
\[\ker^-_b(Z
)= \bigcup_{0 \leq r < b} \ker(\rho^{\scst Z}_{0r})
=\big\langle
\big\{[z_j]+[z_i]\colon (z_j, b_j, z_i)\in \TMT(Z)\mbox{ and }
b_j<b\big\}\big\rangle\]
and $\ker^+_b(Z)
= \ker(\rho^{\scst Z}_{0b})
=\big\langle\big\{
[z_j]+[z_i]\colon (z_j, b_j, z_i)\in \TMT(Z)\mbox{ and }
b_j\leq b\big\}\big\rangle$.
\\
We also define 
$\ker_0^-(Z)=0$ and $
\ker_0^+(Z) =  \big\langle \big\{ 
[z_j] + [z_i]
\mid (z_j, 0, z_i)\in \TMT(Z) 
\big\} \big\rangle $. Finally,
\[
\ker^-_\infty(Z) = \bigcup_{0\leq r} 
\ker(\rho^{\scst Z}_{0r})
\;\;\mbox{ and }\;\;
\ker^+_\infty(Z) = \Ho_0(\VR_0(Z))\,.
\]

We consider now another set of $n$ points 
$X=\{x_1,x_2,\dots,x_n\}$ 
together with a symmetric non-negative function
 $\dX\colon X\times X\rightarrow \R$,
and a bijection $\fbullet\colon X\to Z$ mapping $x_i$ to $z_i$ for all $i\in\{1,2,\dots,n\}$.
Such a bijection induces an isomorphism $f_0\colon\Ho_0(\VR_0(X))\to \Ho_0(\VR_0(Z))$. 
In particular, 
we have that 
$f_0(\ker^-_0(X))=0$, and also
\begin{align*}
  f_0 \big(\ker^-_a(X)\big) = \big\{
[f_0(x_j)]+[f_0(x_i)]\colon (x_j, a_j, x_i)\in \TMT(X)\mbox{ and }
a_j< a
\big\} \mbox{ for all $a>0$, and}
\\
f_0 \big(\ker^+_a(X)\big) = \big\{
[f_0(x_j)]+[f_0(x_i)]\colon (x_j, a_j, x_i)\in \TMT(X)\mbox{ and }
a_j\leq a
\big\}  \mbox{ for all $a\in\R$}.
\end{align*}
The {\it induced block function} 
$\cMXZ\colon \R\times\R
\rightarrow \N$ is
defined, for all
$a,b \in \R$,   by
\begin{equation*}
\cMXZ(a,b) = 
\dim\bigg(
\dfrac{
f_0(\ker_{a}^+(X)) \cap \ker_{b}^+(Z)
}{
f_0(\ker_{a}^-(X)) \cap \ker_{b}^+(
Z) +f_0(\ker_{a}^+(X)) \cap \ker_{b}^-(Z)
}
\bigg)
\end{equation*}
Intuitively,
$\cMXZ(a,b)$ is the amount of connected components of $\Ho_0(\VR_0(X))$ that die at $a$ and are sent by $f_0$ to connected components of $\Ho_0(\VR_0(Z))$ that die at $b$.
It is well-defined because 
$f_0(\ker_{a}^{\pm}(X))\subseteq \Ho_0(\VR_0(Z))$ and 
$\ker_{b}^{-}(Z)\subseteq \ker_{b}^{+}(Z)\subseteq \Ho_0(\VR_0(Z))$.
This definition is an adaptation, to our setting,  of the one given  in~\cite{matchings} whose worst-case complexity is $O(n^3)$ as stated in \cite{tdqual}.
Unlike in~\cite{matchings}, it is possible that $\cMXZ(a,b)\neq 0$ for 
a pair $(a,b)$ with $a<b$,
and we cannot guarantee the 
 very good properties presented in~\cite{matchings}  such as, for example, linearity
as, in general, $f_\bullet$ does not induce a persistence morphism $f\colon \PH_0(X)\rightarrow \PH_0(Z)$.

However,
there is a simple workaround.
Given $\delta\in\R$, we define $\Xdelta$ to be
the set $X$ together with the symmetric non-negative function $d_{\Xdelta}:X\times X\to\R$ given by $\dX^{\delta} (x,x')=\dX(x,x')+\delta$.
Now, for
$\delta \in \R$ big enough, $\fbullet^{\delta}:\Xdelta\to Z$ is a \emph{non-expansive map}, meaning that
$d_{\Xdelta}(x,x')\geq \dZ(\fbullet^{\delta}(x),\fbullet^{\delta}(x'))$ for all $x,x'\in \Xdelta$. 
Taking advantage of such property,  we derive the following result, the proof of which can be found in Appendix~\ref{prooflem:mf}.

\begin{lemma}\label{lem:mf}
  If $\delta\in\R$ is big enough, then 
  $f^{\delta}
  \colon \PH_0(\Xdelta)\to \PH_0(Z)$ is a persistence morphism and we have that
\[
  f^{\delta}
   \simeq 
    \mbox{$\bigg(
        \bigoplus_{ b \in \R} \bigoplus_{a\geq b} 
        \bigoplus_{ i\in \{1,2,\dots,\cM_{f^{\delta}}^0(a,b)\}} 
        (\kappa^{a} \rightarrow \kappa^{b})\bigg) 
       \oplus \big(\kinfty\rightarrow \kinfty\big)$.}
        \]
\end{lemma}

\subsection{Induced Matching Distance}

Given barcodes $(S, m)$ and $(S', m')$, a \emph{matching} is a bijection
$
 \sigma\colon\Rep\, (S,m) \rightarrow \Rep\, (S', m')
 $. 
The following is a key property satisfied by the block function $\cMXZ$, that we use to define the induced matching $\sigma_f^0\colon\Rep B(X)\to \Rep B(Z)$.
Its proof is provided in Appendix~\ref{prooflem:matching}.
 
\begin{lemma}\label{lem:matching}
For all $a,b\in\R$, the block function
$\cMXZ$ 
satisfies that 
\[\mbox{$\sum_{b' \in \R} \cM^0_f(a,b') = \mX(a)$ and $\sum_{a' \in \R} \cMXZ(a',b) = \mZ(b)$.}\]
 \end{lemma}

We now introduce the novel concept of \emph{matching distance} induced by bijections between point sets endowed with symmetric non-negative functions. It is worth mentioning paper~\cite{torras} where the related concept of \emph{matching diagrams} induced by set mappings between finite metric spaces,  was defined. 
Fixed $q\in \N$,  the {\it induced matching distance} is defined as:
\[
\dM^q(\B(X),\B(Z))=
\Bigg(
\sum_{\substack{
(a,\ell) \in \Rep B(X) \\
\sigma_f^0((a,\ell))=(b,\ell')
}}
|a-b|^q\Bigg)^{1/q}
=
\Bigg(\sum_{a,b \in \R}
\cMXZ (a,b) \cdot 
|a-b|^q
\Bigg)
^{1/q}
\,.
\]
Observe that  the
$q$-Wasserstein distance \cite{cohen2010lipschitz}  is always bounded by the induced matching distance. 
The paper's main result states that the induced matching distance is continuous. The proof of this result can be consulted in Appendix~\ref{proofth:continuous}.

\begin{theorem}\label{th:continuous}
Fixed $\veps>0$, 
   consider a bijection $f_\bullet\colon X\rightarrow Z$ such that $|\dX(x,y)-\dZ(f_\bullet(x), f_\bullet(y))| < \veps$ for all $x,y\in X$.
    Then, 
    $\dM^q(B(X), B(Z)) \leq (n
     - 1)^{1/q} \cdot \veps$.
\end{theorem}

We conclude this section by illustrating why, in certain cases, it is preferable to compare barcodes using the induced matching distance rather than the $q$-Wasserstein distance. 

\begin{figure}[ht!]
  \includegraphics[width=0.45\textwidth]{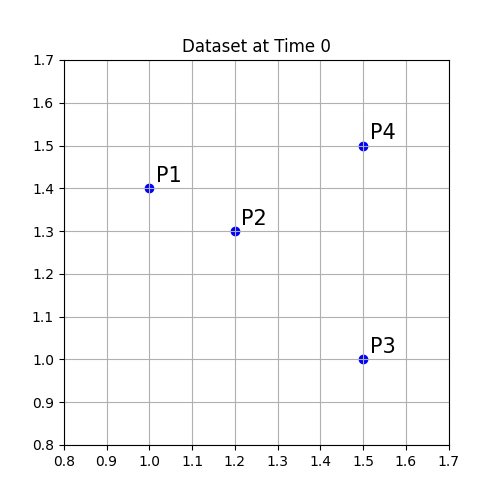}
   \includegraphics[width=0.45\textwidth]{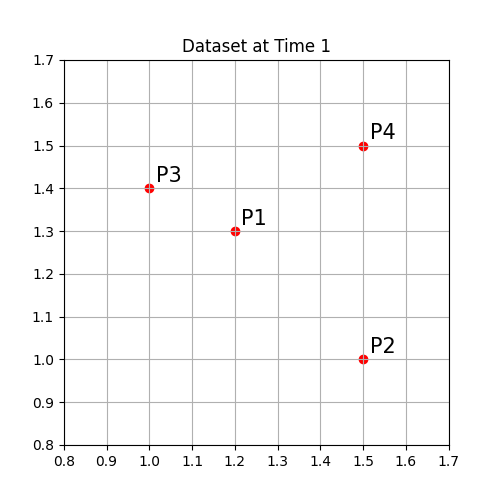}
    \caption{
  Point clouds $X_0$ (left) and $X_1$ (right) with the same points but different labels.
   }
    \label{fig:datasets}
\end{figure}

Consider the point clouds $X_0$ and $X_1$  pictured in Fig.~\ref{fig:datasets}. 
Observe that they have the same shape, and therefore, their barcodes are equal, not reflecting that  the points 
in the clouds have changed 
their labels. 
We aim to preserve the information conveyed by labels since points may represent the states of distinct agents, and the labels identify the agents.

\begin{figure}[ht!]
 \includegraphics[width=0.45\textwidth]{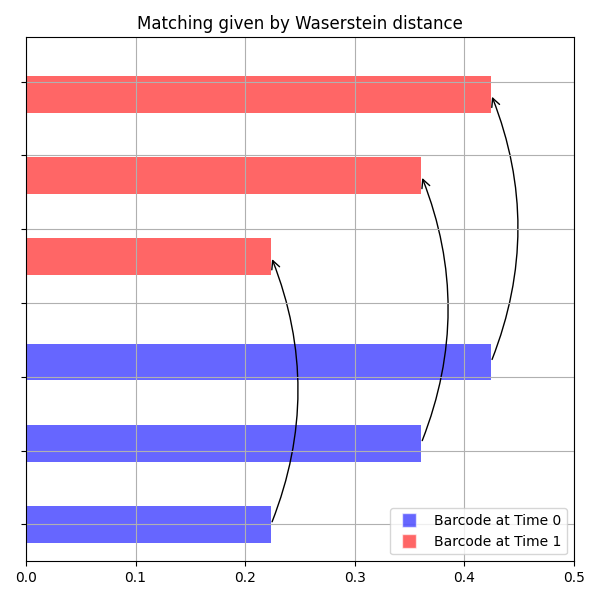}
  \includegraphics[width=0.45\textwidth]{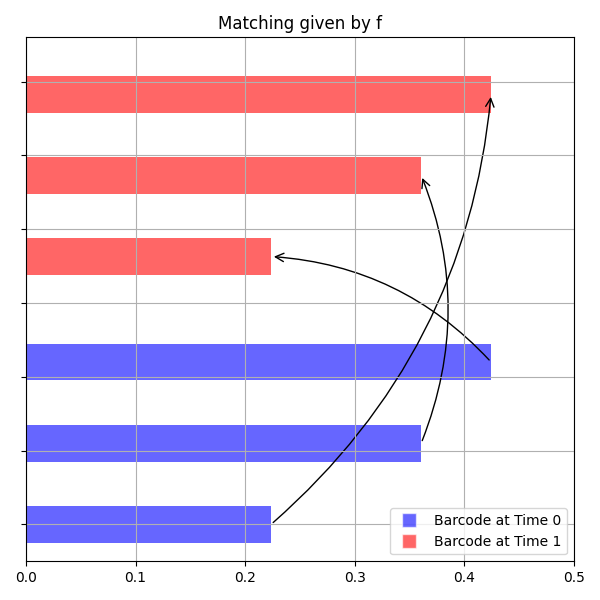}
    \caption{
The matching that produces the $q$-Wasserstein distance (left) and the matching that produces the induced matching distance (right) between $\B(X_0)$ and $\B(X_1)$.
 }
    \label{fig:matchings}
\end{figure}

In Fig.~\ref{fig:matchings} (left), we can 
 see the matching that produces
a $q$-Wasserstein distance of 0.  
In Fig.~\ref{fig:matchings} (right),
 we can see the matching induced by the bijection $X_0\to X_1$  resulting from the label changes, produces 
a nonzero induced matching distance, which is more coherent.

\section{Robot Fleet Navigation Analysis via the Induced Matching Distance}

As an application, we aim to compare, via the induced matching distance,
three local navigation algorithms or \emph{behaviors} for robots: HL, ORCA and SF (see Appendix~\ref{sec:behaviors} for more details).
Simulated robots use these algorithms to reach their targets, while we record their trajectories, composed of 
poses and orientations $(x, y, \alpha)\in 
\mathbb{R}^2 \times [0,2\pi)$.

\begin{figure}[ht!]
 \centering
  \includegraphics[width=0.32\textwidth]{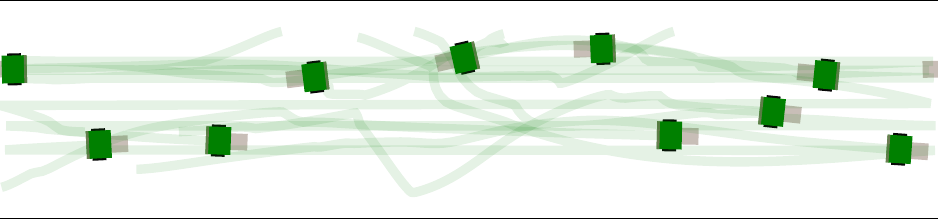}
  \includegraphics[width=0.32\textwidth]{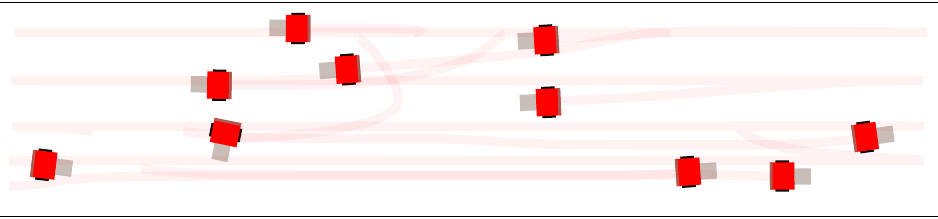}
  \includegraphics[width=0.32\textwidth]{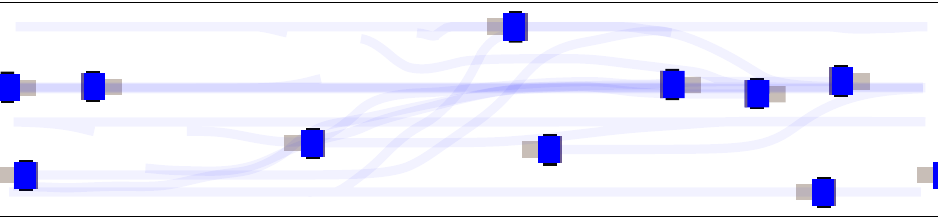}
  \caption{
  The first 10 seconds of the trajectories followed by the robots in the corridor scenario are depicted,  corresponding to the behaviors  SF (left), ORCA (center), and HL (right).
  }
  \label{fig:corridor}
\end{figure}

For our experiments, we use the Navground social navigation simulator \cite{navground}. 
Specifically, we consider simulations in a scenario that represents a corridor with a length of \SI{5}{\meter} and a width of \SI{3.5}{\meter}, having both ends 
connected. In each simulation, 10 autonomous robots apply the same behavior to navigate the corridor, with 5 robots moving to the left and 5 moving to the right. 
These robots mimic smart wheelchairs with differential drive kinematics, radius \SI{0.4}{\meter},  target speed \SI{1.2}{\meter\per\second}, and maximal speed \SI{1.66}{\meter\per\second}. 
At the beginning of the simulations, the robots are randomly distributed in the corridor and attempt to follow their assigned direction while avoiding collisions.
We collect the pose and orientation of each robot every \SI{0.1}{\second} for \SI{90}{\second}, resulting in a set of 10 3-dimensional time series, each containing 900 values. Fig.~\ref{fig:corridor} shows the trajectories followed by the 10 robots for 3 different simulations.

In each simulation, the $i$-th robot is associated with a 3-dimensional time series $a^i$ representing its trajectory, being
$a^i= \big\{a^i_t\big\}_{t=1}^{900}$
 and $a^i_t= (x^i_t,y^i_t,\alpha^i_t)$. Now, for $t=1,2,\ldots,850$, we obtain the set $Z_t = \big\{z^i_t\big\}_{i=1}^{10}$, where $z^i_t=\{a^i_t,a^i_{t+10},\ldots,a^i_{t+50}\} \subset a^i$ is a time window.
We use dynamic time warping (see Appendix~\ref{sec:DTW}) as the symmetric non-negative function to 
build a sequence of Vietoris-Rips filtrations $\big\{\VR_0(Z_t)\big\}_{t=1}^{850}$. Now, for $t=1,2,\ldots,800$, we  compute the induced block $\cM_{f^t}^0$ from the bijection $f^t_{\bullet}: Z_t\to Z_{t+50}$, resulting in a sequence of induced matching distances $\big\{d_{f_0^t}^1(B(Z_t),B(Z_{t+50}))\big\}_{t=1}^{800}$ that is called the \emph{induced matching signal}.

The induced matching signal provides an intuition on the stability of the trajectories.
Ideally, the robots tend to organize themselves and end up forming lanes in the corridor, moving forward indefinitely without ever correcting their trajectories or changing their speed. This is called a stable state. When a simulation has stabilized, the set $Z_t = \{z^1_t,z^2_t,\ldots,z^{10}_t\}$ may change as $t$ changes, but the distances between its elements given by dynamic time warping remain similar, so the Vietoris-Rips filtrations are also similar and the induced matching distances are low. A decreasing induced matching signal indicates that the corresponding simulation is stabilizing, with the robots becoming better organized.

\begin{figure}[ht!]
\centering
\includegraphics[width=0.75\linewidth]{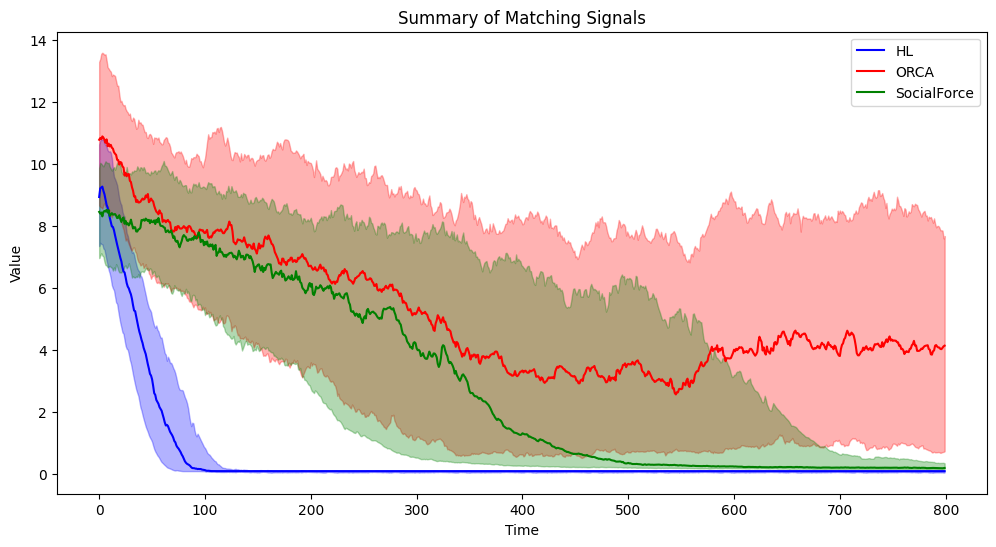}
\caption{Induced matching signals for the 600 simulations in the corridor scenario. Bold lines show the median signal for each behavior, and shaded bands represent the interquartile ranges.}
\label{fig:corridor-matching-signals}
\end{figure}

We have performed 200 simulations for each behavior, getting 600 induced matching signals summarized in Fig.~\ref{fig:corridor-matching-signals}. At first glance, the signals vary significantly depending on the behavior. HL and SF tend to stabilize the simulations in less than \SI{90}{\second}, although at different speeds, while ORCA generally fails to stabilize 
within that time. 

\begin{figure}[ht!]
\centering
\includegraphics[width=1\linewidth]{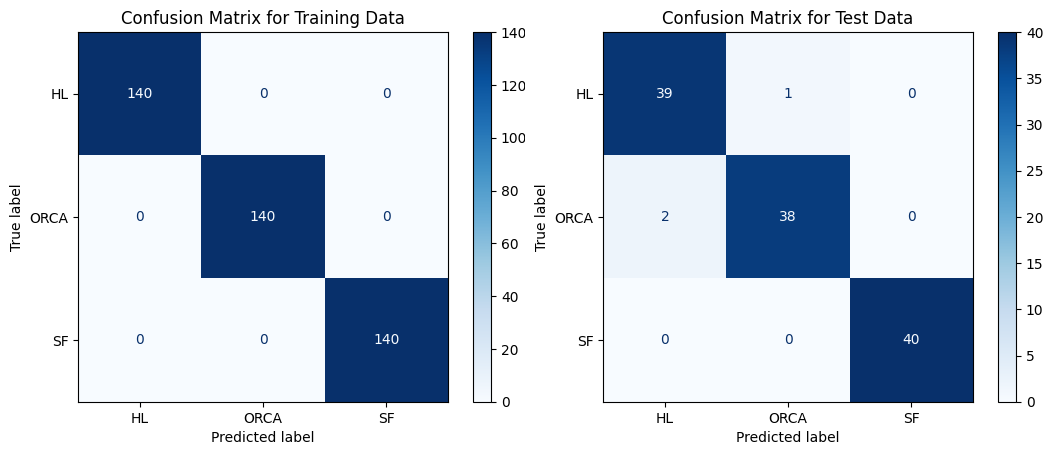}
\caption{The confusion matrices for the ResNet classifier on the induced matching signals.}
\label{fig:corridor-cm}
\end{figure}

To show, in a quantitative way, that the three behaviors can be well distinguished using our topology-based signal, we build a time series classifier based on neural networks.
Following best practices  \cite{ismail2019deep}, we trained a ResNet model \cite{wang2017time} using Python \cite[Chapter~8]{cerqueira2024deep}. 
We split the set of 600 signals into 420 for training, 60 for validation, and 120 for testing. After 100 training epochs, we got a ResNet model whose confusion matrices are shown in Fig.~\ref{fig:corridor-cm}. The model demonstrates perfect performance on the training set and achieves an accuracy and a macro-average F1-score of $97.5\%$ on the test set, indicating that the induced matching signals effectively distinguish between the different navigation behaviors.

\section{Conclusions}

We introduced a topology-based tool for comparing discrete structures represented by symmetric non-negative functions and applied it to monitor the dynamic behavior of trajectory groups.
The study of its stability and extension to higher dimensions remains future work.

\noindent
{\bf Code availability:}
The source code for the experiments  presented in this paper can be found in {\bf \href{https://github.com/Cimagroup/induced-matching-distance-navground}{https://github.com/Cimagroup/induced-matching-distance-navground}.}

\bibliographystyle{plain} 
\bibliography{main}

\appendix

\section{Proof of Lemma~\ref{lemma:decomposition-modules}}\label{prooflemma:decomposition-modules}

First, let 
    $V$ be the persistence module given by 
   \[
    V_r=\Bigg(\oplus_{b \in \R}  \dfrac{\rho_{0r}(\ker^+_b(Z))}{\rho_{0r}(\ker^-_b(Z))}\Bigg) 
    \oplus 
    \Bigg( \dfrac{\rho_{0r}(\ker^+_\infty(Z))}{\rho_{0r}(\ker^-_\infty(Z))}\Bigg)
    \,.
    \]
    for $r \in \R$, with structure maps $\overline{\rho}_{rs}\colon V_r\rightarrow V_s$, for $r<s$, induced by $\rho_{rs}$ on all quotient summands.
    \\
    Now, for all $z_j \in Z$
   with death value 
   $b_j$,
   let $(z_j,b_j,z_i)\in\TMT(Z)$. 
   Then, 
   $\rho_{0r}([z_j]+[z_i]) \in \PH_0(Z)_r$ is sent
    to the 
    quotient 
    class 
    of 
    $[z_j]+[z_i]$ in the summand indexed by $b_j$, while $\rho_{0r}([z_1])$ is sent to the quotient class of $[z_1]$ in the last summand of $V_r$.
    Clearly, this defines linear maps $\psi_r\colon\PH_0(Z)_r \rightarrow V_r$ for all $r \in \R$, 
  that are isomorphisms for all $r\in\R$, and
    such 
    that, for $r<s$, the commutativity $\psi_s \circ \rho_{rs} = \overline{\rho}_{rs}\circ \phi_r$ holds.
            Altogether $\psi\colon \PH_0(Z)\rightarrow V$ is a persistence 
    isomorphism.
    \\
  In addition, there is a persistence isomorphism 
$\phi:     V\to
    \left(\bigoplus_{b \in \R} \bigoplus_{\ell \in 
    \{1,2,\dots,\mZ(b)\}} \kb\right) 
    \oplus 
    \kinfty$.
    \\
    Then,
     we concatenate  $\psi$ with
     $\phi$
       to conclude the proof.

\section{Proof of Lemma~\ref{lem:mf}}\label{prooflem:mf}

If $\delta$ is big enough, then $f^{\delta}_{\bullet}\colon \Xdelta\to Z$ is non-expansive.
In such case, there is an induced graph morphism $\VR_r(X^\delta) \rightarrow \VR_r(Z)$ for all $r\in\R$ inducing a persistence morphism
$f^{\delta}
   \colon \PH_0(\Xdelta)\to \PH_0(Z)$.
Then, we 
 adapt previous results for
persistence morphisms as follows.
First, by 
Theorem~5.3 from~\cite{Jacquard2023}, we have that $f^{\delta}_0$ is isomorphic to:
\[
   \bigg(
        \bigoplus_{ b \in \R} \bigoplus_{a\geq b} 
        \bigoplus_{s=1}^{r^b_a} 
        (\ka \rightarrow \kb)\bigg)     
       \oplus
\bigg(\bigoplus_{c\in\R} 
\bigoplus_{s=1}^{d^-_c}
(0 \rightarrow \kc)\bigg)
\oplus
\bigg(\bigoplus_{c\in\R} \bigoplus_{s=1}^{d^+_c}
(\kc \rightarrow 0)\bigg) 
       \oplus\big(\kinfty\rightarrow \kinfty\big)
        \]
       for some  $r^b_a,d^\pm_b\in\R$, since 
 there are no nested intervals (because 
all intervals of the barcodes $\B(\Xdelta)$ and 
$\B(Z)$ start at $0$). 
Since the bijection $f_\bullet\colon X\to Z$  induces an isomorphism
$f^{\delta}_0:\Ho_0(\VR(\Xdelta))\to \Ho_0(\VR(Z))$,
then $d^\pm_c=0$ for all $c\in\R$, concluding that:
\[f^{\delta}\simeq
   \bigg(
        \bigoplus_{ b \in \R} \bigoplus_{a\geq b} 
        \bigoplus_{s=1}^{r^b_a} 
        (\ka \rightarrow \kb)\bigg)     
              \oplus\big(\kinfty\rightarrow \kinfty\big)
        \]
Now, by Theorem 5.1  from \cite{matchings},
we have that $\cM_{f^{\delta}}^0$ is linear, which implies that
\[
\cM^0_{f^\delta}(c,d) =
\Bigg(\sum_{b \in \R} \sum_{a\geq b} 
r^b_a \cdot 
\cM^0_{\ka\to \kb}(c,d)\Bigg) + \cM^0_{\kappa^\infty\to \kappa^\infty}(c,d)
\]
for all $c,d \in \R$.
By the definition of the block function, $\cM^0_{\ka\to \kb}(c,d)=1$ for $c=a$ and $d=b$, and is zero otherwise. 
Altogether, given $c\geq d$, 
it follows that 
$r^d_c=\cM^0_{f^\delta}(c,d)$, 
and the result holds.

\section{Proof of Lemma~\ref{lem:matching}}\label{prooflem:matching}

First, 
 consider $\delta\in\R$ big enough so that $f^\delta_\bullet:\Xdelta\to Z$ is non-expansive. Then,    $f^{\delta}:
\PH(X^\delta)\to \PH(Z)
  $ is a persistence morphism by Lemma~\ref{lem:mf}. 
  Now,  a consequence of 
  Lemma~\ref{lemma:decomposition-modules} and Lemma~\ref{lem:mf} is that 
  \[\mbox{$\sum_{b' \in \R} \cM_{f^{\delta}}^0(a,b') = \mXdelta(a)$ and $\sum_{a' \in \R} \cM_{f^{\delta}}^0(a',b) = \mZ(b)$.}\]
 Then,  we can always define a matching 
$ \sigma_1\colon\Rep B(\Xdelta)\to \Rep B(Z)$.

Second,  
from the equality $\ker^\pm_a(X)=\ker^\pm_{a+\delta}(X^\delta)$ for all $a \in \R$, we have
\[\begin{tikzcd}
    \phi \colon \SX \ar[r, "\simeq"] & 
     S^{\scst \Xdelta}
     \end{tikzcd}\mbox{ is such that }
     \phi(a)=a+\delta,\;\mbox{  and }\;
     \mX(a)=m^{\scst \Xdelta}(a+\delta),
     \mbox{ for all }a\in\SX\,.
     \] 
  In addition,  
  $\cMXZ(a,b) = \cM^0_{f^{\delta}}(a+\delta, b)$, for all $a,b \in \R$.
Then, we have another matching $\sigma_2\colon\Rep B(X)\to \Rep B(\Xdelta)$. 

We then have the matching 
$ \sigma_1\circ \sigma_2=\sigma_f^0\colon \Rep B(X)\to \Rep B(Z)$, concluding the proof.

\section{Proof of Theorem~\ref{th:continuous}}\label{proofth:continuous}

Fixed $\veps>0$, 
we consider two sets of $n$ points, 
$X$ and $Z$, with symmetric non-negative functions, $\dX$ and $\dZ$, and a bijection $f_\bullet\colon X\rightarrow Z$ such that $|\dX(x,y)-\dZ(f_\bullet(x), f_\bullet(y))| < \veps$ for all$x,y\in X$.
 First, let us prove that 
$\cM^0_{f}(a,b)=0$ for all $a,b \in \R$ such that $|a-b|>\veps$.

Since $|a-b|<\veps$, either $a+\veps < b$ or $b+\veps < a$. Consider the first case. 
    Given $[x_i]+[x_j]\in \ker^+_a(X)$, it follows that there exist a sequence 
          $y_1, y_2, \ldots, y_m$ with $y_1=x_i$, $y_m=x_j$ and $(y_i,y_{i+1})
    \in \VR_a(X)$
    for all $1 \leq i < m$. 
    Then,     
    $[f_\bullet(x_i)] + [f_\bullet(x_j)] \in \ker^+_{a+\veps}(Z)$ since, by hypotheses, $\dXp(f_\bullet(y_i), f_\bullet(y_{i+1}))< a+\veps$ for all $1 \leq i < m$.
    Now, since 
    $f_0(\ker^+_a(X)) \subseteq \ker^+_{a+\veps}(Z)\subseteq \ker^-_{b}(Z)$, it follows that $\cMXZ(a,b)=0$. 
    The case when $b+\veps < a$ is analogous, where one uses that $f_\bullet$ is a bijection to show that $\ker^+_{b}(Z) \subseteq f_0(\ker^-_{a}(X))$, which also implies the result.
We conclude that, if $\cMXZ(a,b)\neq 0$ then $|a-b|<\veps$.    

 Now, using Lemma~\ref{lem:matching}, we have that
        $\sum_{b \in \R} \cMXZ(a,b)$ $ = \mX(a)$.
     And using the fact that $
           n-1 = \sum_{a \in \R}\mX(a)$, we obtain: 
    \[
        \dM^q(\B(X),\B(Z)) 
        \leq 
        \Bigg(\sum_{a,b \in \R}\cMXZ (a,b)\cdot    
        \veps^q\Bigg)^{1/q}
        = 
        \Bigg(\sum_{a \in \R}
        \mX(a)\cdot     
        \veps^q\Bigg)^{1/q}
        = \Big( (n-1) \cdot \veps^q \Big)^{1/q}
    \]
    and the result holds.

\section{Robot social navigation simulation}\label{sec:behaviors}

In all three navigation algorithms considered, 
individual robots take control actions operating on the local environment around them, which contains the state of their neighbors and the position of surrounding obstacles. 
There is no explicit coordination between the robots: possible collective behaviors, like the formation of ordered lanes of flow, emerge from the interaction of the individual behavior without an explicit design or centralized coordination.

\begin{itemize}

    \item The {\it social force (SF)} behavior \cite{SocialForceBehavior} describes pedestrian dynamics using a physics-inspired model based on fictitious forces. In this model, the pedestrians are driven by a force that directs them towards their desired goal, while other forces repel them from obstacles and other pedestrians to avoid collisions. This method effectively captures 
    social interactions of pedestrians in crowded environments, simulating realistic and natural trajectories. 

    \item The {\it optimal reciprocal collision avoidance (ORCA)} behavior \cite{ORCAbehavior}  is a computationally efficient local navigation algorithm used in robotics and crowd simulation. ORCA belongs to the family of algorithms based on the Reciprocal Velocity Obstacle~\cite{van2008reciprocal}, where agents first compute a set of safe velocities, and then select the nearest safe velocity to their target velocity, and where pairs of agents on collision path, share half of the responsibility to avoid the collisions.
   ORCA is designed for agents to avoid collisions 
     cooperatively and efficiently by adjusting their velocities based on the positions and velocities of other agents in their immediate environment.

    \item The {\it human-like (HL)} behavior \cite{HLikeJerome} is a bio-inspired, computationally light, local navigation algorithm for robotics, that adapts a heuristic model for pedestrian motion. It addresses engineering aspects like trajectory effectiveness and scalability, as well as societal aspects by producing human-friendly, predictable trajectories. Like SF, equipping 
    a robot with HL promotes its friendliness and predictability around people, but contrary to SF, HL searches for safe control actions: first selecting the velocity direction coming nearest to the target, and then fixing the speed in order to prevent collisions.
    In other words, this behavior is inspired by how pedestrians move.

\end{itemize}

\section{Dynamic Time Warping}~\label{sec:DTW}

Instead of comparing the pose 
at each moment to understand the behavior of two robots, we can look at their entire trajectories by comparing the full paths of poses 
in a sliding window, giving a better understanding of how the robots behave. 
This way, considering all the steps in a trajectory in a time window, we get two three-dimensional time series, one for the poses and the other for the twists:

\[
\left[\begin{array}{c}
\mbox{step } 1\\
\mbox{step } 2\\
\vdots
\\
\mbox{step } n
\end{array} 
\right] \Rightarrow 
\left[\begin{array}{c}
(x,y,\alpha)_1\\
(x,y,\alpha)_2\\
\vdots
\\
(x,y,\alpha)_n
\end{array} 
\right]\;\;
\mbox{ and }\;\;
\left[\begin{array}{c}
(v_x,v_y,v_{\alpha})_1\\
(v_x,v_y,v_{\alpha})_2\\
\vdots
\\
(v_x,v_y,v_{\alpha})_n
\end{array} 
\right]
\]
We could compare the difference between trajectories using similarity functions, such as the Euclidean distance or the dynamic time warping, between the corresponding time series.

Dynamic time warping (DTW) is a technique widely recognized in the speech recognition community. It enables a non-linear alignment between two signals by minimizing the distance between them. 
The DTW concept was introduced in \cite{DTW-classic}
for spoken word recognition and was used in \cite{DTW-second}
to find patterns in time series. In both cases, the time series were 1-dimensional. 
DTW was later adapted for multi-dimensional signals in different ways as noticed in \cite{adaptative-dtw} where an adaptive approach is considered.

The classical technique of DTW uses a dynamic programming approach to align two given time series. Specifically, to compute the dynamic time warping of the time series $S$ and $T$  of 
length $m$ and $n$, respectively, namely $S=(s_1,s_2,\dots,s_m)$ and $T=(t_1,t_2,\dots,t_n)$, we first form a $n\times m$ grid $G$.
Each grid point $(i,j)$ for $i\in\{1,2,\dots,n\}$ and 
$j\in\{1,2,\dots,m\} $
corresponds to a distance measure
 between the elements $s_i$ and $t_j$, that is
 $G(i,j)=d(s_i,s_j)$
 for a given distance $d$. 
 A warping
path $W$  aligns the elements of $S$ and $T$ such that the total distance between their elements is
minimized.
More formally, $W$ is a path of ${\cal A}(S,T)$ that minimizes the following value:
\[
DTW(S,T)=
\min_{\pi\in {\cal A}(s,t)}
\sum_{(i,j)\in\pi} d(s_i,s_j)
\]
 were ${\cal A}(S,T)$ is the set of paths that align the elements of $S$ and $T$. 
 Specifically, a sequence $\pi=(\pi_1,\pi_2,\dots,\pi_k)$ is a path of ${\cal A}(S,T)$ if:
 \begin{itemize}
          \item for all $\alpha\in\{1,2,\dots,k\}$ we have that  $\pi_{\alpha}=(i_{\alpha},j_{\alpha})$
     with $i_{\alpha}\in \{1,2,\dots,n\}$ and
     $j_{\alpha}\in \{1,2,\dots,m\}$.
     \item $\pi_1=(1,1)$ and $\pi_k=(n,m)$.
     \item The sequence $\pi$ is 1-Lipschitz  and monotonically increasing, that is,
     for $\pi_{\alpha}=(i_{\alpha},j_{\alpha})$ and 
     $\pi_{\alpha+1}=(i_{\alpha+1},j_{\alpha+1})$ we have that 
     $i_{\alpha+1}-i_{\alpha}\leq 1$  and 
     $  j_{\alpha+1}-j_{\alpha}\leq 1$.
 \end{itemize}

To illustrate why DTW might be more suitable for comparing trajectories than Euclidean distance, 
think in two trajectories,
with one trajectory displaced relative to the other. 
Observe that DTW aligns points in the time series that are locally similar, unlike the alignment used to compute the Euclidean distance. 
On the other hand, unlike Euclidean distance, DTW is not a metric because two different trajectories can have a distance of zero, and the triangle inequality does not always hold.

\end{document}